\documentclass[12pt, reqno]{amsart}
\usepackage{amsmath, amsthm, amscd, amsfonts, amssymb, graphicx, color}
\usepackage{setspace}
\usepackage{mathrsfs}
\usepackage{multicol}
\usepackage[bookmarksnumbered, colorlinks, plainpages]{hyperref}
\hypersetup{colorlinks=true,linkcolor=red, anchorcolor=green, citecolor=cyan, urlcolor=red, filecolor=magenta, pdftoolbar=true}

\newtheorem{theorem}{Theorem}[section]
\newtheorem{lemma}[theorem]{Lemma}

\newtheorem{corollary}[theorem]{Corollary}
\theoremstyle{definition}

\theoremstyle{remark}

\numberwithin{equation}{section}

\newcommand{\NN}{\mathbb{N}}

\newcommand{\CC}{\mathbb {C}}

\begin{document}
\setcounter{page}{1}
\title[Cyclic and supercyclic  weighted composition operators on the Fock space]{Cyclic and supercyclic weighted composition operators on the Fock space }
\author[Tesfa  Mengestie]{Tesfa  Mengestie }
\address{ Mathematics Section \\
Western Norway University of Applied Sciences\\
Klingenbergvegen 8, N-5414 Stord, Norway}
\email{Tesfa.Mengestie@hvl.no}


 \begin{abstract}
 We study the cyclic and supercyclic dynamical properties of weighted composition operators on the  Fock space $\mathcal{F}_2$.  A complete characterization of cyclicity  which depends on the derivative of the symbol for  the composition operator and  zeros of the weight function  is  provided. It is  further shown that  the space fails to support supercyclic  weighted composition operators. As a consequence, we also noticed that the space supports no cyclic multiplication operator.
\end{abstract}
\maketitle
\section{Introduction}
For a pair of entire functions   $(u, \psi)$  on the complex plane  $\CC$,  the  induced  weighted composition operator $W_{(u,\psi)}$ maps $ f$ to $u f(\psi).$  If  $u=1$, then    $W_{(u,\psi)}$ is just  the composition  map $C_\psi: f\mapsto f(\psi)$. On the other hand,  if $\psi$ is the identity map, then  $W_{(u,\psi)}$ reduces to the multiplication operator $M_u:  f\mapsto  uf$. Thus,  $W_{(u,\psi)}$ generalizes the two operators and can be   also  written as a  product $ M_u C_\psi$. The theory of weighted composition operators traces back to the sixties in the  work of
Forelli \cite{For} where it was shown that  the isometries in the Hardy spaces $ H^p$ whenever $1 < p <\infty , p \neq 2 $ are weighted composition operators. De Leeuw \cite{Hof} later    showed that the  same  holds true on the space $H^1$ as well. Since then the operator  became a natural object of study  and  its investigations  has  rapidly evolved in function related operator theory. A number of researchers have studied the operator over various settings with the aim to express its spectral, topological and dynamical properties in terms of the function theoretic properties of the inducing pairs of   symbols $(u,\psi)$. See for example \cite{CZ,Tle,TM4,UK} and the references therein.

In this note we study the  cyclic and supercyclic dynamical properties of the operators  on the classical Fock space $\mathcal{F}_2$. Recall that $\mathcal{F}_2$  consists of  square integrable analytic functions in $\CC$ with respect to the Gaussian measure $d\mu(z)= \frac{1}{\pi}e^{-|z|^2}  dA(z)$ where $dA$ is the Lebesgue measure in  the complex plane $\CC$.  It is a reproducing kernel Hilbert space with inner product
\begin{align*}
\langle f, g\rangle=   \int_{\CC} f(z)\overline{g(z)}d\mu(z),
\end{align*}norm  $\|f\|_2:=\sqrt{ \langle f, f\rangle} $, and kernel function $K_w(z)= e^{\langle z, w\rangle}$ .

 A great deal of the studies on weighted composition operators has  been devoted to characterizing  boundedness and compactness spectral properties over various functional spaces.  On Fock spaces, these properties have been well understood  and described for example in \cite{Tle, TM4,UK} and  expressed in different conditions among which following \cite{Tle},   $W_{(u,\psi)}$ is bounded on $\mathcal{F}_2$ if and only if
 \begin{align}
 \label{bounded}
 \sup_{z\in \CC} |u(z)|e^{\frac{1}{2}(|\psi(z)|^2-|z|^2)} <\infty.
 \end{align}Furthermore, it was   proved that \eqref{bounded} implies $
 \psi(z)= az+b, |a|\leq 1$
 and whenever $|a|=1$, the multiplier function  has the special  form
 \begin{align}u(z)=u(0)K_{-\overline{a}b}(z)
 \end{align} for all $z$ in $\CC$.
  Similarly, compactness has been described by the additional  condition that
    $|a|<1.$
    An interesting  feature here is that there exists an interplay between  $u$ and $\psi$ such that $W_{(u,\psi)}= M_u C_\psi$ is bounded (compact) on $\mathcal{F}_2$ while both
    $C_{\psi}$  and $u$ fail to be. As pointed earlier, the purpose of this note is to study the effect of this interplay  on the dynamical structures of  $W_{(u,\psi)}$ on $\mathcal{F}_2$. Recall that  a bounded linear operator $T$ on a separable  Banach space $\mathcal{H}$ is said to be cyclic  if there exists a vector $f$ in $\mathcal{H}$ for which the   span of the orbit  \begin{align*}\text{Orb}(T,f)=\big\{ f,\  Tf,\  T^2f,\  T^3 f, \  ...\big\}\end{align*}  is dense in $\mathcal{H}$. Such a vector is called cyclic  for $T$.  The operator is hypercyclic if  the   orbit it self  is dense.  $T$ is supercyclic with vector $f$ if  the projective orbit
 \begin{align*}\CC.\text{ Orb}(T,f)= \big\{ \lambda T^nf, \ \ \lambda \in \CC, \ n= 0, 1, 2, ...\big\}\end{align*}
 is dense.
 These dynamical properties of  $T$ depend on the behaviour of its iterates $T^n= T \circ T\circ T \circ...\circ T-n \ \text{times}$. For detailed backgrounds, one may  consult  the  monographs
 \cite{BM,Grosse}.

  It is worth noting that identifying cyclic and  hypercyclic operators have been  a subject of high interest partly because they play central rolls in the study of other operators. More specifically, it is known
that  every bounded linear operator
on an infinite dimensional  complex separable  Hilbert space is the sum of two hypercyclic
operators \cite[p. 50]{BM}. Interestingly, this result holds true with the summands being cyclic operators \cite{WU}. In \cite{GK,TMW} it was  reported that there exists no supercyclic (and hence hypercyclic) composition operator on Fock spaces. On the other hand,  the orbit of any vector $f$ under  $W_{(u,\psi)}$  has elements of the form
\begin{align}
\label{interplay}
W_{(u,\psi)}^n f= f(\psi^n) \prod_{j=0}^{n-1} u(\psi^j)
\end{align} for all nonnegative integers $n$ and $\psi^0$ is the identity map. This shows that the product of weighted composition operators is another weighted composition operator with symbol  $(U_n, \psi^n)$ where
\begin{align*}
 U_n= \prod_{j=0}^{n-1} u(\psi^j).
\end{align*}
 The formula in \eqref{interplay} further  displays a kind of interplay between the  functions $\psi $ and $u$ and generates interest to ask  whether the interplay results in  hypercyclic,   supercyclic or cyclic  weighted composition operators on the Fock space. In fact,  we can make the following simple  observation right away.
  Boundedness of $W_{(u,\psi)}$ implies that  $\psi(z)= az+b, |a|\leq 1$. If $a\neq 1$, then $\psi$ fixes the point $z_0= b/(1-a)$ in $\CC$.  A  computation with adjoint property gives
\begin{align*}
 W_{(u,\psi)}^* K_{z_0}= \overline{u({z_0})} K_{\psi({z_0})}=\overline{u({z_0})} K_{z_0}.
\end{align*}
Since $ K_{z_0}$ is a nonzero vector, it follows that $\overline{u({z_0})}$ is an eigenvalue for  the adjoint operator  $W_{(u,\psi)}^*$ which contradicts the general fact that the adjoint of a hypercyclic operator can not have eigenvalue (see Proposition 1.17 in \cite{BM}). This shows there exists no hypercyclic weighted  composition operator on $\mathcal{F}_2$  in this case. Here the fixed point behaviour of the composition operator decisively determined the absence of hypercyclicity for the weighted operator.   Later,  we will in fact show that the same conclusion applies when $a= 1$ as well.  Composition operators on Fock spaces are among the groups  of non-hypercyclic operators.  It turns out that the weighted composition operators
exhibit the exact same hypercyclic  phenomena as the unweighted ones. The relation in \eqref{interplay} is weak enough  to make any orbit as big as the whole space.

 Given the above observation on the absence of hypercyclic weighted composition operator $W_{(u,\psi)}$, the natural question is  what happens to the weaker   cyclicity  and supercyclicity properties. Clearly,  hypercyclicity is a much stronger property than cyclicity and every hypercylic operator is cyclic. Thus, hypercyclic operators enjoy richer operator-theoretic properties than cyclic ones.   One interesting difference between the two properties is that if an operator has a hypercyclic vector $f$, then it has a dense subset of hypercylic vectors because every element in its orbit is also hypercyclic.  This fails to hold for  cyclic operators: see \cite{SHP} for an example.  On the other hand, supercyclicity is a property which is intermediate between the two.
 \section{main results}
 We now state our first main result on the cyclic behavior of $W_{(u,\psi)}$.
\begin{theorem}\label{thm1}
Let $u$ and $\psi(z)=az+b, \ |a|\leq 1$ be analytic maps on $\CC$  such that $W_{(u,\psi)}$  is bounded on  $\mathcal{F}_2$. Then $W_{(u,\psi)}$
  is cyclic on $\mathcal{F}_2$ if and only
     \begin{enumerate}
  \item  $u$ fails to vanish on $\CC$ and
  \item $a^k\neq a$ for all positive integer $ k\geq2$.
  \end{enumerate}
  \end{theorem}

 The cyclic property of $W_{(u,\psi)}$ depends  on  the  size of the power of   the derivative of the symbol
   $\psi(z)= az+b$ and the existence of zeros of the  symbol $u$. Part b) is the condition required for the cyclicity of the composition operator \cite{GK,TMW}. It means that the absence of zeros for the multiplier function makes it possible for  the weighted operator  to inherent the same  condition. Notice that the cyclicity  condition clearly  restricts  $\psi$  to be a nonconstant function.
 \begin{proof}
 For the sufficiency, we may exhibit that  $K_z$ is a cyclic vector for  $W_{(u,\psi)}$. Since condition ii) implies that $0\neq a\neq1$, observe that using \eqref{interplay} we may rewrite
 \begin{align*}
  W_{(u,\psi)}^m K_z(w)=  K_z\bigg( a^mw+ \frac{1-a^m}{1-a}b\bigg) \prod_{j=0}^{m-1} u(\psi^j(w)) \quad \quad \quad \quad \quad \quad \quad \quad \\
  = e^{\big( a^mw+ \frac{1-a^m}{1-a}b\big)\overline{z}} \prod_{j=0}^{m-1} u(\psi^j(w)) = \sigma_m K_{\overline{a}^m z}(w)  \prod_{j=0}^{m-1} u(\psi^j(w)) ,
 \end{align*} where
 \begin{align*}
 \sigma_m=  e^{\big(  \frac{1-a^m}{1-a}\big)b\overline{ z}}
 \end{align*}  is independent of the point $w$.
  Now,  fix a vector  $f\in \mathcal{F}_2$ that is orthogonal to the orbit of $K_z$ under $W_{(u,\psi)}$. We claim that $f$  should be  the zero function. To  see this, note that for  every $m$
  \begin{align}
  \label{vanish}
  0= \langle f, \  W_{(u,\psi)}^m K_z \rangle= \overline{\sigma_m}\bigg\langle f,\  \prod_{j=0}^{m-1} u(\psi^j) K_{\overline{a}^m z}\bigg\rangle \quad \quad \quad \quad \quad \quad \quad \quad\nonumber \\
  =\overline{\sigma_m}\bigg\langle \prod_{j=0}^{m-1} \overline{u(\psi^j)} f, K_{\overline{a}^m z}\bigg\rangle= \overline{\sigma_m}\prod_{j=0}^{m-1} \overline{u(\psi^j(\overline{a}^m z))} f(\overline{a}^m z).
  \end{align}  By assumption, since $u$ posses  no zero so does the quantity
  \begin{align*}
  \prod_{j=0}^{m-1} \overline{u(\psi^j(\overline{a}^m z))}.
  \end{align*} If  $|a|<1 $, then  $\sigma_m$ is different from zero for every $m$, and  since $\overline{a}^m z \to 0$
 \begin{align*}
 \sigma_m(z)= e^{ \big(\frac{1-a^m}{1-a}\big)b\overline{z}} \to  e^{ \frac{b\overline{z}}{1-a}}
 \end{align*} as $m\to \infty$ for each fixed point $z$.  Thus, the relation in \eqref{vanish} holds only if  $f$ vanishes  on the null  sequence $\overline{a}^m z$.  It follows that $f$ is the zero function   as asserted.

  For the case  $|a|=1$,  and $a^m\neq a$ for each $m\geq 2$, we may argue as follows to arrive at the same conclusion, namely that $f$ is still the zero function.  For each fixed $z$ and all $m\geq 0$, we have
 \begin{align}
 \label{zero}
f(\overline{a}^m z)=0.
 \end{align} Set in particular $z= 1$  and  consider the Taylor series expansion of $f$ at it;
 \begin{align}
 \label{tylor}
 f(w)= \sum_{n=0}^\infty a_n (w-1)^n.
 \end{align}  We need to show that $a_n= 0$ for all $n$. Applying \eqref{tylor} and \eqref{zero} for $m= 0$, gives that
 $a_0= 0$.  Using this and similarity considering   $m= 1, 2, 3, 4, 5, ...$  we get an infinite system of linear equations $VA= 0$ where
  \begin{align*}
  \begin{split}
  V=
\begin{pmatrix}
            1 & \overline{a}-1& (\overline{a}-1)^2& (\overline{a}-1)^3 &... \\
           1 & \overline{a}^2-1& (\overline{a}^2-1)^2& (\overline{a}^2-1)^3 &...\\
           1 & \overline{a}^3-1& (\overline{a}^3-1)^2& (\overline{a}^3-1)^3 &...\\
           . & . &. &.& ...\\
            . & . &. &.& . ..\\
             . & . &. &. &...
         \end{pmatrix}
         \text{ and }
         A=
         \begin{pmatrix}
          a_1 \\
          a_2 \\
          a_3\\
          .\\
          .\\
          .\end{pmatrix}
          \end{split}
           \end{align*}
Observe that  this is a matrix with the terms of a geometric progression in each row, and restricting  to  finite $n\times n$  actually gives the well-known  Vandermonde matrix
\begin{align*}
V_n= \begin{pmatrix}
            1 & \overline{a}-1& (\overline{a}-1)^2& (\overline{a}-1)^3 &... &(\overline{a}-1)^{n-1} \\
           1 & \overline{a}^2-1& (\overline{a}^2-1)^2& (\overline{a}^2-1)^3 &...& (\overline{a}^2-1)^{n-1}\\
           1 & \overline{a}^3-1& (\overline{a}^3-1)^2& (\overline{a}^3-1)^3 &...& (\overline{a}^3-1)^{n-1}\\
           . & . &. &.& . & .\\
            . & . &. &.& . & .\\
             . & . &. &. &. & .\\
             1 & \overline{a}^{n-1}-1& (\overline{a}^{n-1}-1)^2 & (\overline{a}^{n-1}-1)^3 &...&(\overline{a}^{n-1}-1)^{n-1}
                      \end{pmatrix}
\end{align*} with
determinant
 \begin{align}
 \label{vander}
 \det V_n= \prod_{1\leq i<j\leq n} (\overline{a}^j-\overline{a}^i).
 \end{align}
Since $a^m \neq a$ for all $m>1$, we have that $\overline{a}^m-1\neq 0$ for all $m>0$. In fact for two different positive integers $p$ and $q, \ \overline{a}^p\neq \overline{a}^q$. This shows that \eqref{vander} is non-zero.  Thus, the  Vandermonde matrix is invertible and  the solution of the above mentioned system of linear equations is unique for any dimension $n$ of the
system. But  $n$ is arbitrary above, therefore, the infinite system  $VA= 0$ holds only if the system has a trivial solution   $A=0$  as claimed and completes the proof of the sufficiency.

We next proof the necessity of the conditions, and assume on the contrary  that $u$
vanishes at the point $z_0$ in $\CC$. Then,   for any possible cyclic vector $f$ again
 \begin{align*}
W_{(u,\psi)}^n f(z_0)= \prod_{j=0}^{n-1} u\big(\psi^j\big(z_0\big)\big). f\big(\psi^n\big(z_0\big)\big)= u(z_0) \prod_{j=1}^{n-1} u\big(\psi^j\big(z_0\big)\big). f\big(\psi^n\big(z_0\big)\big)= 0
\end{align*}
which shows that all the functions in the orbit  vanish at the point
$z_0$. This extends to all functions $g$ in the closed linear span of the orbit of $f$ under $W_{(u,\psi)}$ which is not the case.

For part b), we consider first the case when $a= 1$ and hence $\psi^j(z)= z+jb$ for all $j\geq 0$.  Then for any possible cyclic vector $f$
\begin{align*}
 W_{(u,\psi)}^m f(z)= f(z+mb) u(0)^m\prod_{j=0}^{m-1} K_b(z+jb)= f(z+mb) u(0)^m e^{ h_m(z) },
\end{align*} where
\begin{align*}
h_m(z):= -\overline{b}\sum_{j=0}^{m-1}(z+jb)= -\overline{b}mz-\frac{|b|^2}{2} m(m-1)=-\frac{|b|^2}{2} m(m-1) K_{-mb}(z).
\end{align*}
First observe that if $b=0$, then the relations above give
\begin{align}
W_{(u,\psi)}^m f= f u(0)^m
\end{align} asserting that all the elements in the orbit are scalar multiplies of the cyclic vector $f$. Any  vector $g$ in $ \mathcal{F}_2$ orthogonal to $f$ is also orthogonal to the closed linear span of its orbit. Thus,  it suffices to show that there exists such a non-zero $g$.  To this end, since $S=\{f\}$ is a closed subspace of  $\mathcal{F}_2$, the projection operator $P: H\to S$ is continuous. Then for  $f_1\in \mathcal{F}_2, f_1\neq f$, the function \begin{align*}g= f_1-Pf_1= f_1-f\in S^\perp\end{align*}  and orthogonal to the linear span of the orbit.

Next, we  assume  $b\neq0$, and consider $g\in \mathcal{F}_2 $   such that
\begin{align}
\label{contra}
\langle g, \  W_{(u,\psi)}^m f \rangle= 0,
\end{align} for all $m$.  Our proof will be completed if we manage to construct a nonzero function $g\in \mathcal{F}_2$  satisfying condition \eqref{contra}.  Taking the condition further
\begin{align*}
0=\langle g,  W_{(u,\psi)}^m f \rangle= \overline{u(0)}^m e^{-\frac{|b|^2}{2} m(m-1)}\langle g \overline{f(\psi^m)},K_{-mb} \rangle\\
= \overline{u(0)}^m e^{-\frac{|b|^2}{2} m(m-1)} g(-mb) \overline{f(\psi^m(-mb))} \\
=\overline{u(0)}^m e^{-\frac{|b|^2}{2} m(m-1)} g(-mb) \overline{f(0)}
\end{align*}  which holds only if $g$ vanishes on the sequence $\{-bm: m\in \NN \}$.   A good candidate is to set  $g$ to be the product
\begin{align*}
g(z)= z\prod_{m=1}^\infty \bigg(1+\frac{z}{bm}\bigg) e^{-\frac{z}{bm}}.
\end{align*} Accordingly, it suffices to show that the product above converges and belongs to the space $\mathcal{F}_2$. But it is know that the function $\sin z$ has  the product expansion
\begin{align*}
\sin z= z\prod_{m=1}^\infty \bigg(1-\frac{z}{m}\bigg) e^{\frac{z}{m}}
\end{align*}  from which we may rewrite  $g$ as
\begin{align*}
g(z)= \frac{-b}{\pi} \sin\Big(\frac{-\pi  z }{b}\Big),
\end{align*} and observe that  $g$ belongs to the space with norm  $\|g\|_2\leq |b|/\pi$.

 It remains to show the case when $\psi(z)= az+b, |a|=1$ and $a\neq1$. In this case

$ \psi$  fixes the point $z_0= \frac{b}{1-a}$.  A simple modification of the arguments in the proof of Lemma~3 of \cite{GG}  gives that the adjoint operator $W_{(u,\psi)}^*$  has the following  set of eigenvalues
\begin{align}
\label{eigen}
\Big\{\overline{u(z_0)}, \ \overline{au(z_0)}, \ \overline{a^2u(z_0)}, \  \overline{a^3u(z_0)}, ..., \ \overline{a^mu(z_0)}\Big\},
\end{align} where $m$ is the smallest positive integer such that $m\geq 2$ and $a^m= a$.  Note that such an $m$ exists by our assumption.  On the other hand,  $u(z_0)\neq 0$ and hence the  set in \eqref{eigen} has at least two elements.  By Proposition 2.7 of \cite{SHP}, the adjoint of a cyclic operator can not have multiple eigenvalues from which our assertion follows and completes the proof of the first result.
 \end{proof}
Unlike cyclicity, our next result shows that the space $\mathcal{F}_2$ fails to support supercyclic weighted composition operators.
\begin{theorem}\label{thm2}
Let $u$ and $\psi$ be  analytic  maps on $\CC$  such that $W_{(u,\psi)}$  is bounded on  $\mathcal{F}_2$. Then $W_{(u,\psi)}$ can not be  supercyclic    on $\mathcal{F}_2$.
\end{theorem}
 The result extends earlier work obtained for the unweighted composition operators on the space \cite{TMW}. The interplay between  the additional multiplier function $u$ and $\psi$ has failed to make any projective orbit dense in  $\mathcal{F}_2$.
It is natural to compare this result with the  analogues space over the unit disc where it  was reported that there exists hypercyclic and hence supercyclic weighted composition operators, see for example \cite{YY}. It seems these properties are less attainable  on functional spaces defined on unbounded domains than bounded ones.
\begin{proof}
Since $W_{(u,\psi)}$  is bounded, we set $\psi(z)
=az+b, \ |a|\leq 1$.  As  supercyclicity implies cyclicity, the proof for  the case when  $a^m= a$ for some $m>1$  follows from  Theorem~\ref{thm1}.  Therefore, we plan to prove  the case for
$a^m\neq a $ for all $m\geq 2$.
 Assume first that  $|a|<1$. Then   $W_{(u,\psi)}$ is a compact operator. For this we rely on some  spectral properties of supercyclic operators.  The space  $\mathcal{F}_2$ is an infinite dimensional complex Banach  space, and if a compact operator $W_{(u,\psi)}$ is supercyclic on  $\mathcal{F}_2$, then    its  spectrum $ \sigma(W_{(u,\psi)})$  contains only the zero element: see \cite [p.29]{BM}. Thus, it suffices to show that $ \sigma(W_{(u,\psi)})$ contains  at least two  elements.  Since  the point $z_0= b/(1-a)$  is fixed by $\psi$, we simply change the kernel functions and   repeat the arguments in the proof of Theorem~1 in \cite{GG} to  deduce
\begin{align}
\label{spectrum}
\sigma(W_{(u,\psi)})=\big\{0, \ u(z_0), \ au(z_0), \ a^2u(z_0), \ a^3u(z_0),... \big\}.
\end{align}
Observe that if  $u(z_0)= 0$, then the operator is not supercyclic, If not, $W_{(u,\psi)}$ becomes cyclic and contradicts Theorem~\ref{thm1}.  On the other hand, if $a= 0$, then  $\psi(z)= b$ and  the relation in  \eqref{interplay} implies
 \begin{align*}
 W_{(u,\psi)}^n f(z)= f(b) u(b)^n
 \end{align*} This means that the orbit of $f$ contains only constant functions, and hence $W_{(u,\psi)}$ is not supercyclic here again.  Therefore,  both $a$ and $u(z_0) $ are nonzero and
 the set in \eqref{spectrum} contains infinitely many elements.

It remains to show the case for $|a|=1$ and $a^m \neq a$ for all $m>1$.  This is rather immediate as  $u(z)=  u(0) K_{b}(z) $, by  Theorem 3.2 of \cite{Tle}, $W_{(u,\psi)}$ is a constant multiple of a unitary operator and hence normal.  Then,  our conclusion follows from a result of Hilden and Wallen \cite{HW}.
\end{proof}
The following corollary is now immediate.
 \begin{corollary}
 \begin{enumerate}
 \item Let $C_\psi$ be bounded on $\mathcal{F}_2$, that is $\psi(z)= az+b, |a|\leq 1$ and $b= 0$ whenever $ |a|= 1.$ Then $C_\psi$ is cyclic if and only if $a^k\neq a$ for  each positive integer $k\geq 2$. In this case, the reproducing kernel $K_z$ is  a cyclic vector.
     \item $C_\psi$  is  not supercyclic on $\mathcal{F}_2$.
     \item There exists no cyclic  multiplication operator  on $\mathcal{F}_2$.

     \end{enumerate}
   \end{corollary}


His research interests are in functional analysis, complex analysis, and   operator theory. More interested at the interface of these three branches of mathematics.

\end{document}